\documentclass[12pt,draft]{article}
\usepackage{amssymb, amsmath, amsthm,  graphicx}
\usepackage[cp1251]{inputenc}
\usepackage[english]{babel}
\usepackage{cite, multirow, enumerate, array}
\newcounter{parag}
\newcommand{\sect}[1]
{\refstepcounter{parag}
\begin{center} {\bf \theparag. #1} \end{center}}

\newtheorem{cor}{Corollary}[parag]
\newtheorem{theorem}{Theorem}[parag]
\newtheorem{lemma}{Lemma}[parag]
\newtheorem{prop}{Proposition}[parag]

\theoremstyle{definition}
\newtheorem*{prf}{Proof}

\begin{document}

\begin{center} {\bf \large On element orders in covers of finite \\ simple groups of Lie type}\\[0.5cm]
{\bf Mariya A. Grechkoseeva}\footnote{Supported by the RFBR (12-01-31221 and 12-01-90006).}
\end{center}

\sect{Introduction}

Given a finite group $G$, we say that a finite group $H$ is a {\it cover} of $G$ if $G$ is a homomorphic image of $H$. Denoting by $\omega(A)$ the set of element
orders of a group $A$ and referring to it as the {\it spectrum} of $A$, we say that $G$ is {\it recognizable (by spectrum) among covers} if $\omega(G)\neq\omega(H)$ for any proper cover $H$ of $G$. If $G$ has a nontrivial
normal soluble subgroup, then $G$ is not recognizable among covers, and in fact, there are infinitely many nonisomorphic covers $H$ of $G$ with $\omega(H)=\omega(G)$ \cite[Lemma 1]{98Maz.t}.

The present paper concerns recognizability of finite nonabelian simple groups, and we denote these groups according to the `Atlas of finite groups' \cite{85Atlas}. The following simple groups are known to be recognizable by spectrum among covers: the sporadic groups (see \cite{98MazShi}), the alternating groups \cite{ 99ZavMaz.t}, the Ree and Suzuki groups \cite{92Shi, 93BrShi, 99DenShi}, the exceptional groups of types $G_2$ and $E_8$ \cite{13VasSt.t,10Kon}, the~groups $L_2(q)$ \cite{94BrShi, 99ZavMaz.t},
the groups $L_3(q)$ \cite{04Zav, 04Zav1.t}, the groups $L_n(q)$ with $n\geqslant 5$ \cite{08Zav1.t}, and
the~groups $U_n(q)$, with $n\geqslant 5$, other than $U_5(2)$ \cite{11Gr}. Also it is worth noting that $L_4(q)$ and $U_4(q)$ are recognizable among covers as long as $q$ is even or prime \cite{08MazChe.t, 08Zav1.t, 11Gr}, and $U_3(q)$ is not recognizable among covers if and only if $q$ is a Mersenne prime with $q^2-q+1$ being a prime too \cite{06Zav.t}. Lastly, $L_4(13^{24})$, $U_5(2)$ and $^3D_4(2)$ are not recognizable among covers \cite{08Zav2, 11Gr, 13Maz.t}.

Our purpose is not to give an exhaustive answer to the question what finite simple groups of Lie type are recognizable by spectrum among covers and which are not but
rather to show that groups of sufficiently large Lie rank have this property of recognizability. Thus we do not touch the very interesting case of the groups $L_4(q)$ and $U_4(q)$ but focus our attention on groups not mentioned above.
Recall that $O_{2n+1}(q)\simeq S_{2n}(q)$ provided that $n=2$ or $q$ is even.

\begin{theorem}\label{t:main}
If $S$ is one of the simple groups $^3D_4(q)$ with $q>2$, $F_4(q)$, $E_6(q)$, $^2E_6(q)$, $E_7(q)$, $S_{2n}(q)$ with $n\geqslant 2$, $O_{2n+1}(q)$ with $n\geqslant 3$ and $q$ odd, and $O_{2n}^\pm(q)$ with $n\geqslant 4$, then $S$ is recognizable by spectrum among covers.
\end{theorem}

The theorem and before-mentioned results yield the following corollary.

\begin{cor}
If $S$ is a finite nonabelian simple group other than $L_4(q)$, $U_3(q)$, $U_4(q)$, $U_5(2)$, and $^3D_4(2)$, then $S$ is recognizable by spectrum among covers.
\end{cor}

Our strategy of proof is standard and based on the following two observations (in these observations, as well as throughout the paper, all modules and representations are finite-dimensional). First, it is not hard to see that a group $G$ is recognizable by spectrum among covers as long as $\omega(G)\neq \omega(V\leftthreetimes G)$ for any nontrivial $G$-module $V$ over a field of positive characteristic (cf. Lemma \ref{l:reduction} below). Second, given  a finite group of Lie type, it is natural to distinguish representations in the defining characteristic from other representations. Thus for every simple group $S$ satisfying the hypothesis of the theorem and every nontrivial $S$-module $V$, we need to prove that $\omega(S)\neq \omega(V\leftthreetimes S)$ provided that the corresponding representation (i) is equicharacteristic, and (ii) is cross-characteristic (cf. Problems 17.73 and 17.74 in \cite{10Kou}).

\begin{prop}\label{p:equi}
Let  $q$ be a power of a prime $p$ and let  $S$ be one of the simple groups  $^3D_4(q)$ with $q>2$, $F_4(q)$, $E_6(q)$, $^2E_6(q)$, $E_7(q)$, $S_{2n}(q)$ with $n\geqslant 2$, $O_{2n+1}(q)$ with $n\geqslant 3$ and $q$ odd, and $O_{2n}^\pm(q)$ with $n\geqslant 4$. If $V$ is a nonzero $S$-module over a field of characteristic $p$ and $H$ is a~natural semidirect product of $V$ and $S$, then $\omega(S)\neq \omega(H)$.
\end{prop}

Under condition (ii), the problem has already been settled for symplectic and orthogonal groups of dimension at least six \cite{11Gr}. So we deal with the remaining groups.

\begin{prop}\label{p:cross}
Let  $S$  be one of the simple groups  $^3D_4(q)$, $F_4(q)$, $E_6(q)$, $^2E_6(q)$, $E_7(q)$, and $S_{4}(q)$. If $V$ is a nonzero $S$-module over a field of characteristic $r$, where $r$ is coprime to $q$, and $H$ is a natural semidirect product of $V$ and $S$, then $\omega(S)\neq \omega(H)$.
\end{prop}

Proposition \ref{p:equi} and \cite{13VasSt.t} give an affirmative answer to  \cite[Problem 17.73]{10Kou} in all cases except for the cases (a) and (f). In the case (f) the answer is negative \cite{13Maz.t} and the case (a)
concerns $U_4(q)$. Proposition \ref{p:cross}, \cite[Theorem 1]{11Gr} and \cite[Lemma 12]{08Zav1.t}  give an affirmative answer to \cite[Problem 17.74]{10Kou}.

We conclude this paper with a result comparing the spectrum of a proper cover of $O_{2n+1}(q)$ with the spectrum of $S_{2n}(q)$. This result is motivated by
the following assertion: if $L=S_{2n}(q)$, where $n\geqslant 4$, and $G$ is a finite group such that $\omega(G)=\omega(L)$, then $G$ has a unique
nonabelian composition factor and if this factor is a group of Lie type in the same characteristic as $L$, then it is isomorphic either to $L$ or to one of the groups $O_{2n+1}(q)$ and $O_{2n}^-(q)$ \cite[Theorem 3]{09VasGrMaz.t}.
It is well-known that the spectra of $S_{2n}(q)$ and $O_{2n+1}(q)$ are very close. This property of $S_{2n}(q)$ and $O_{2n+1}(q)$ makes the case when the factor is isomorphic to $O_{2n+1}(q)$ more difficult to eliminate. Somewhat paradoxically, the same property allows us to apply our results on covers to this situation.

\begin{prop}\label{p:bncn}
Let $L=S_{2n}(q)$ and $S=O_{2n+1}(q)$, where $n\geqslant 3$ and $q$ is odd. If $H$ is a proper cover of $S$, then $\omega(H)\not\subseteq \omega(L)$.
In particular, if $G$ is a finite group such that $\omega(G)=\omega(L)$ and $G$ has a composition factor isomorphic to $S$, then $S\leqslant G\leqslant \operatorname{Aut}(S)$.
\end{prop}

\sect{Preliminaries}

Our notation for classical groups and groups of Lie type is based on that in \cite{98GorLySol} with the~following exception: $\Phi(q)$ always stands for the adjoint version,
while the universal version is denoted by $\Phi(q)_u$. The same rule applies to the symbol $^d\Phi(q)$. Also
sometimes we substitute $\Phi(q)$ and $^2\Phi(q)$ by $\Phi^+(q)$ and $\Phi^-(q)$ respectively .

\begin{lemma}[Zsigmondy \cite{Zs}]\label{l:zsig}
Let $q\geqslant 2$ and $n\geqslant 3$ be integers, with $(q,n)\neq (2,6)$. There exists
a prime $r$ such that $r$ divides $q^n-1$ but does not divide $q^i-1$ for $i<n$.
\end{lemma}

With notation of Lemma \ref{l:zsig}, we call a prime $r$ a \textit{primitive prime divisor} of $q^n-1$ and denote it by $r_n(q)$.

By $(n,m)$ we denote the greatest common divisor of two positive integers  $n$ and $m$.
Given a positive integer $n$, we write $\pi(n)$ for the set of all prime divisors of $n$.
Also, if $r$ is another positive integer, then $n_r$ denotes the $\pi(r)$-part of $n$ (that is, the largest divisor $m$ of $n$ such that
$\pi(m)\subseteq \pi(r)$), and $n_{r'}$ denotes the ratio $n/n_r$.

If $G$ is a finite group and $r$ is a prime, we refer to the highest power of $r$ that lies in $\omega(G)$ as the $r$-{\it exponent} of $G$.
If $G$ is a group of Lie type over a field of positive characteristic $p$, then it is natural and instructive to consider $\omega(G)$ as a disjoint union of
the following three subsets: the subset of numbers that are powers of $p$, the subset $\omega_{p'}(G)$ of numbers that are coprime to $p$, and the subset $\omega_{mix}(G)$ of the rest, so called `mixed', orders.

If $G$ is a finite group, then $\pi(G)=\pi(|G|)$. The {\it prime graph } of $G$ is a graph whose vertex set is $\pi(G)$ and in which two different
vertices are adjacent if and only if their product lies in $\omega(G)$. 
The following property of the simple groups of Lie type is well-known but rarely appears as a~separate assertion, so we include the corresponding lemma in this paper.

\begin{lemma}\label{l:pg}
Suppose that $S$ is a finite simple group of Lie type. Then for every $r\in\pi(S)$, there is $t\in\pi(S)$ such that $r\neq t$ and $rt\not\in\omega(S)$.
\end{lemma}

\begin{prf}
If $S$ is one of the simple groups $L_2(q)$, $L_3(q)$, $U_3(q)$, and $S_4(q)$, or $S$ is a simple exceptional group of type other than $E_7$, then
by \cite{81Wil, 89Kon.t}, the prime graph of $S$ is disconnected, and the~assertion follows. For the rest of simple classical groups, the lemma is proved in \cite[Lemma~12]{11Gr}.
If $S=E_7(q)$, then we can take $t$ to be equal to $r_9(q)$ or $r_{18}(q)$ (see, for example,  \cite[Fig.~4]{11VasVd.t}).
\end{prf}

The next lemma is also well-known and in the special case $A=B$, is proved, for example, in \cite[Lemma 12]{04Zav}.

\begin{lemma}\label{l:reduction}
Let $A$ and $B$ be finite groups. The following are equivalent.
\begin{enumerate}
\renewcommand{\labelenumi}{\rm (\theenumi)}

\item $\omega(H)\not\subseteq\omega(B)$ for any proper cover $H$ of $A$;

\item $\omega(H)\not\subseteq\omega(B)$ for any split extension $H=K:A$, where $K$ is a nontrivial elementary abelian group.

\end{enumerate}
\end{lemma}

\begin{prf}
Clearly, (1) implies  (2). Suppose that (2) does not imply  (1), and let $H$ be a group of minimum possible order such that $H$ is a proper cover of $A$ and $\omega(H)\subseteq\omega(B)$.
Let $K$ be a normal subgroup of $H$ with $H/K=A$. Choose $r\in\pi(K)$ and let $R$ be a Sylow $r$-subgroup of $K$. If $N=N_H(R)$, then the Frattini argument yields
$A\simeq H/K=NK/K\simeq N/N\cap K$. Since $N\cap K\neq 1$ and $\omega(N)\subseteq \omega(H)$, we have $H=N$ by the minimality of $H$. Repeating this observation for other prime divisors of $|K|$,
we conclude that $K$ is nilpotent. Denote the~Hall $r'$\nobreakdash-subgroup of $K$ by $C$ and the Frattini subgroup of $R$ by $\Phi(R)$. The factor group of $H/(C\times \Phi(R))$ by $K/(C\times \Phi(R))$ is isomorphic to $A$, and again the minimality of $H$ yields $C=\Phi(R)=1$. Thus
$K$ is an elementary abelian $r$-group, and in particular, the group $A=H/K$ acts on $K$ by conjugation. If $H_1$ is a semidirect product of $K$ and $A$ with respect to this action, then $\omega(H_1)\subseteq \omega(H)\subseteq\omega(B)$
by \cite[Lemma 10]{99ZavMaz.t}. This contradicts (2).
\end{prf}

\begin{lemma}[{\rm \cite[Lemma 1]{97Maz.t}}] \label{l:frob_action} Let $G$ be a finite group, $K$ a normal subgroup of $G$
and let $G/K$ be a Frobenius group with kernel $F$ and cyclic complement $C$. If $(|F|, |K|)=1$ and
$F$ is not contained in $KC_G(K)/K$, then $r|C|\in\omega(G)$ for some prime divisor $r$ of $|K|$.
\end{lemma}

\begin{lemma} \label{l:frob_groups} The group $SL_n(q)$, where $n\geqslant 2$, includes a Frobenius subgroup with kernel of order $q^{n-1}$ and
cyclic complement of order $(q^{n-1}-1)_{(n,q-1)'}$, as well as a Frobenius subgroup with kernel of order  $q^{k}$ and
cyclic complement of order $q^{k}-1$ for every $1\leqslant k<n-1$.
\end{lemma}

\begin{prf}
It is well known that the natural semidirect product of a vector space of order $q^k$
by the~Zinger cycle of order $q^k-1$ in $GL_k(q)$ is a Frobenius group. It is easily seen that this product can be
embedded into a line stabilizer in $GL_{k+1}(q)$. Since $GL_{k+1}(q)$, where $1\leqslant k<n-1$, can be embedded into
$SL_n(q)$, the second assertion follows. Furthermore, $PSL_n(q)$ includes a~Frobenius subgroup $F$ with kernel of order
$q^{n-1}$ and cyclic complement of order $(q^{n-1}-1)/d$, where $d=(n,q-1)$ (see, for example, \cite[Lemma 5]{05VasGr.t}). Denoting $\pi=\pi(d)$ we see that the full preimage of
the Hall $\pi'$-subgroup of $F$ in $SL_n(q)$ is a split extension of a group of order $d$ by the desired Frobenius group of order $q^{n-1}(q^{n-1}-1)_{d'}$.
\end{prf}

\begin{lemma} \label{l:frob_groups2} The simple group $S_{2n}(q)$, where $q$ is even and $n=2^m\geqslant 2$, includes a Frobenius subgroup with cyclic kernel of order $q^n+1$
and cyclic complement of order $2n$.
\end{lemma}

\begin{prf}
The proof is similar to that of Lemma 2.4 in \cite{04VasGr.t}.
\end{prf}

\begin{lemma}\label{l:hh} Let $S$ be a finite simple group of Lie type over a field of characteristic $p$ and let $S$ acts faithfully on a
vector space  $V$ over a field of characteristic $r$, where $r\neq p$. Let $H=V\leftthreetimes S$ be a natural semidirect product of $V$ by $S$.
Suppose that $s$ is a power of $r$ and that $s$ lies in the~spectrum of some proper parabolic subgroup $P$ of $S$. If the unipotent radical of $P$ is abelian or  both $p$ and $r$ are odd or $p=2$ and $r$ is not a Fermat prime or $r=2$ and $p$ is not a Mersenne prime, then $rs\in\omega(H)$.
\end{lemma}

\begin{prf} Let $g\in P$ be an element of order $s$, $U$ a unipotent radical of $P$ and let $M=U\leftthreetimes \langle g\rangle$.
By  \cite[Theorem 2.6.5e]{98GorLySol} it follows that $C_G(U)=Z(U)$. Thus $O_r(M)=1$, and by the Hall--Higman theorem
\cite[Theorem 2.1.1]{56HalHig},  the degree of the minimal polynomial of $g$ on $V$ is equal to $s$. This implies that the coset $Vg$ contains an element of order $rs$.
\end{prf}

\sect{SR-subgroups of finite groups of Lie type}

Let $G$ be a finite universal group of Lie type. By  \cite[Theorem 2.2.6e]{98GorLySol}, there are a simply connected simple algebraic group $\overline{G}$ and a Frobenius endomorphism $\sigma$ of $\overline{G}$ such that
$G=C_{\overline{G}}(\sigma)$. By \cite[Theorem 2.1.6]{98GorLySol}, there is a $\sigma$-stable maximal torus  $\overline{T}$ of $\overline G$. Let $\Phi$ denote the root system of $\overline G$ with respect to
$\overline{T}$ and $\Pi$ denote some  fundamental system. Also we write $\overline X_\alpha$ to denote a root
subgroup of $\overline G$ associated with a root $\alpha$. Let $\Delta\subseteq \Pi$.
Following \cite{07SupZal}, we call a subgroup of $\overline G$ an {\it SR-subgroup} (with respect to $\overline T$) if it is generated by the subgroups $\overline X_{\pm\alpha}$ for  $\alpha\in\Delta$.
A subgroup of the finite group $G$ is an $SR$-subgroup if it is an intersection of an~$SR$-subgroup of $\overline{G}$ with $G$.

\begin{lemma}\label{l:sr}
Let $G$ be a finite universal group of Lie type over a field of characteristic $p$ and let $g\in G$ lie in some  $SR$-subgroup of type $\Psi$. The element $g$ has an nonzero fixed vector in every nonzero $G$-module over a field of characteristic $p$ in the following cases:

\begin{enumerate}[{\rm (a)}]

\item $G=A_n(q)_u$, $n\geqslant 1$, and $\Psi=A_m$, where $m\leqslant n/2$;

\item $G=B_n(q)_u$,  $n\geqslant 3$, $q$ is odd,  and  $\Psi=A_{n-1}$;

\item $G=C_n(q)_u$,  $n\geqslant 2$, $q$ is odd, and $\Psi=C_{n-1}$;

\item $G=D_n^\pm(q)_u$,  $n\geqslant 4$, and $\Psi=A_{n-2}$;

\item $G=E_7(q)_u$ and $\Psi=E_6$;

\item $G=E_6^\pm(q)_u$ and $\Psi=D_5$.
\end{enumerate}

\end{lemma}

\begin{prf}

See Corollaries 1.6, 1.7, 4.3, 5.1, and also a remark after Corollary 1.3 in \cite{07SupZal}.
\end{prf}

We need an easy consequence of Lemma \ref{l:sr} for simple groups. For convenience, we list generators of the centers
of some simple algebraic groups from \cite[Table 1.12.6]{98GorLySol}.
In Table~\ref{tab:z}, the symbol $\overline G$ denotes a simple algebraic group of type $\Phi$, and $h_{i}(t)$ denotes $h_{\alpha_i}(t)$.

\begin{table}[h]
\caption{Generators of $Z(\overline G)$ for some $\overline G$}\label{tab:z}
\begin{center}
\renewcommand{\arraystretch}{1.3}
\begin{tabular}{|l|m{140pt}|l|}

\hline
$\Phi$& Dynkin diagram & generators of $Z(\overline G)$\\
\hline
$B_n$&\includegraphics[draft=false, scale=0.6]{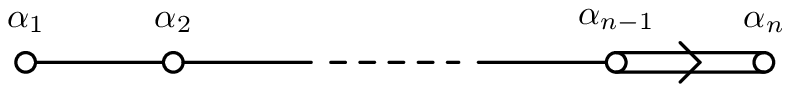}& $h_{n}(-1)$\\
$C_n$&\includegraphics[draft=false, scale=0.6]{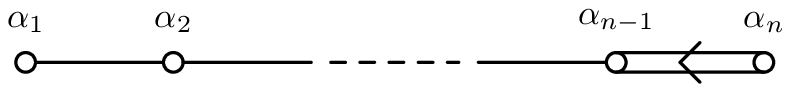}& $h_1(-1)h_3(-1)\dots h_{n}((-1)^n)$\\
$D_{n}$, $n$ even& \multirow{2}{*}{\includegraphics[draft=false, scale=0.6]{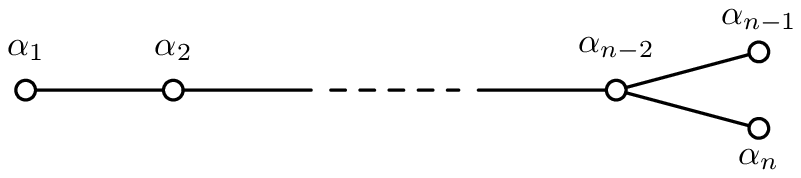}} &$h_{1}(-1)h_{3}(-1)\dots h_{{n-1}}(-1)$,  $h_{{n-1}}(-1)h_{{n}}(-1)$\\
$D_{n}$, $n$ odd && $h_{1}(-1)h_{3}(-1)\dots h_{{n-2}}(-1)h_{{n-1}}(t)h_{{n}}(t)$, $t^2=-1$\\
$E_6$&\includegraphics[draft=false, scale=0.6]{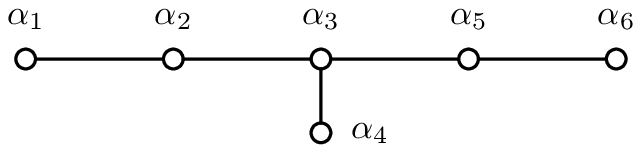}& $h_{1}(t)h_{2}(t^2)h_{{5}}(t)h_{{6}}(t^2)$, $t^3=1$\\
$E_7$&\includegraphics[draft=false, scale=0.6]{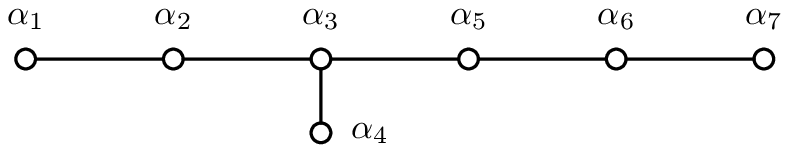}& $h_{4}(-1)h_{5}(-1)h_{{7}}(-1)$\\
\hline
\end{tabular}
\end{center}
\end{table}

\begin{lemma}\label{l:spectrumofc}
Let $S$ be a finite simple groups of Lie type over a field of characteristic $p$, $V$ a nonzero $S$-module over a field of the same characteristic $p$ and $H$ the natural semidirect product of $V$ by $S$.
Then $\omega(H)$ includes the set $p\cdot\omega_{p'}(\Sigma)$, where

\begin{enumerate}[{\rm (a)}]

\item $\Sigma=SL_n(q)$ if $S=O_{2n+1}(q)$, $n\geqslant 3$, $q$ is odd or $S=O_{2n+2}^\pm(q)$, $n\geqslant 3$;

\item $\Sigma=SU_4(q)$ if $S=O_9(q), O_{10}^\pm(q)$,  $q$ is odd;

\item $\Sigma=E_6(q)_u$ if $S=E_7(q)$;

\item $\Sigma=Spin_{10}(q)$ if $S=E_6(q)$;

\item $\Sigma=Spin_8^-(q)$ if $S={}^2E_6(q)$.
\end{enumerate}

\end{lemma}

\begin{prf}
Let $G$ be a finite universal group of Lie type such that $G/Z(G)=S$. We may assume that $V$ is a $G$-module. We prove the lemma by finding an  SR-subgroup $\Gamma$ of $G$ such that
$\Gamma\simeq \Sigma$, $\Gamma\cap Z(G)=1$ and every element of $\Gamma$ has a nonzero fixed vector in $V$.

Since SR-subgroups are a special case of subsystem subgroups, we need some information on the subsystem subgroups from \cite[Section 2.6]{98GorLySol}. Let $\overline G$, $\sigma$, $\overline T$, $\Phi$, and $\Delta$ have the same meaning as those in the beginning of this section. Observe that $Z(G)=Z(\overline G)\cap G$ by  \cite[Proposition 2.5.9]{98GorLySol}.
Let $\tau$ be the isometry of order 2 of the system $E_6$ induced by the symmetry of the~Dynkin diagram if $S={}^2E_6(q)$, and take $\tau=1$ in other cases. Supposing that $\Delta$ generates closed, irreducible and   $\tau$-invariant subsystem $\Psi$ of $\Phi$, we define $\overline\Gamma_\Delta=\langle \overline X_\alpha \mid \alpha\in\Psi\rangle$.
Since $\langle \overline X_\alpha \mid \alpha\in\Psi\rangle=\langle \overline X_\alpha, \overline X_{-\alpha} \mid \alpha\in\Delta\rangle$ by \cite[Theorem 1.12.7d]{98GorLySol}, it follows that $\overline\Gamma_\Delta$ is an SR-subgroup.
Let $W$ and $W_0$ be the Weyl groups of $\Phi$ and $\Psi$ respectively, and let $W_0^*$ be the stabilizer in  $W$ of $\Psi$. For every $w\in W_0^*$ we define an SR-subgroup $\Gamma_{\Delta,w}$ as follows. Let $\overline N=N_{\overline G}(\overline T)$ and choose $n\in\overline N$ mapping to $w$. Then choose $x\in\overline G$ so that $\sigma=(\sigma n)^x$ and set $$\Gamma_{\Delta,w}=C_{\overline \Gamma_\Delta ^x}(\sigma).$$ If $\tau w$ induces an isometry of $\Psi$ of order $d_1$, then  $\Gamma_{\Delta,w}$ is the universal group of type $^{d_1}\Psi(q)$ by \cite[Proposition 2.6.2b,d]{98GorLySol}.

Let $S=O_{2n+1}(q)$, where $n\geqslant 3$ and $q$ is odd. Take $\Delta=\{\alpha_1,\alpha_2,\dots,\alpha_{n-1}\}$ with notation for roots as in Table~\ref{tab:z} and set $\Gamma=\Gamma_{\Delta,1}$. Then every element of $\Gamma$ has a nonzero fixed vector in $V$ by Lemma~\ref{l:sr}. Also it follows by the previous paragraph that $\Gamma=A_{n-1}(q)_u=SL_n(q)$. Lastly, $Z(G)=\langle h_{\alpha_n}(-1)\rangle$ and $G$ is universal, therefore, $\Gamma\cap Z(G)=1$. The groups $O_{2n+2}^\pm(q)$, $n\geqslant 3$, $E_7(q)$, $E_6(q)$ and $^2E_6(q)$ can be handled similarly
with subsystems of types $A_{n-1}$, $E_6$, $D_5$ and $D_4$, respectively, in place of $\Psi$.

Let $S$ be as in (b). Since $O_9(q)$, where $q$ is odd, can be embedded into $O_{10}^\pm(q)$, it suffices to consider the case $S=O_9(q)$. Let $\Delta=\{\alpha_1,\alpha_2,\alpha_{3}\}$ (again with notation of Table~\ref{tab:z}). It is not hard to see that $|W^*_0/W_0|=2$ and some $w\in W_0^*\setminus W_0$ induces an isometry of $\Psi$ of
order two. Taking $\Gamma=\Gamma_{\Delta,w}$ and reasoning as before, we conclude that $\Gamma={}^2A_3(q)_u=SU_4(q)$, every element of $\Gamma$ has a nonzero fixed vector in $V$, and $\Gamma\cap Z(G)=1$. This completes the proof.
\end{prf}

\sect{Proof of Theorem \ref{t:main}}

We begin work by proving Proposition \ref{p:equi}.  Following \cite{01BabSha}, we call a finite group $L$ of Lie type over a field of characteristic $p$ {\it unisingular} if
for every nontrivial abelian $p$-group $A$ with $L$-action, every element of $L$ has a nontrivial fixed point on $A$.

\setcounter{prop}{0}
\setcounter{parag}{1}

\begin{prop}
Let  $q$ be a power of a prime $p$ and let  $S$ be one of the simple groups  $^3D_4(q)$ with $q>2$, $F_4(q)$, $E_6(q)$, $^2E_6(q)$, $E_7(q)$, $S_{2n}(q)$ with $n\geqslant 2$, $O_{2n+1}(q)$ with $n\geqslant 3$ and $q$ odd, and $O_{2n}^\pm(q)$ with $n\geqslant 4$. If $V$ is a nonzero $S$-module over a field of characteristic $p$ and $H$ is a~natural semidirect product of $V$ and $S$, then $\omega(S)\neq \omega(H)$.
\end{prop}

\begin{prf} If the action of $S$ on $V$ is not faithful, then $C_S(V)$ is a nontrivial normal subgroup of $S$, and hence $C_S(V)=S$. This implies that for every prime
$t\in\pi(S)\setminus\{p\}$, the product $pt$ lies in $\omega(H)$, and so $\omega(H)\neq\omega(S)$ by Lemma \ref{l:pg}. Thus we can assume that $S$ acts faithfully on $V$.

Let $S=F_4(q)$. By \cite{84Der}, we know that $q^4-q^2+1\in\omega(S)$ and $p(q^4-q^2+1)\not\in\omega(S)$. Since $S$ is unisingular \cite[Theorem 1.3]{03GurTie}, it follows that $p(q^4-q^2+1)\in\omega(H)\setminus\omega(S)$.

Let $S={}^3D_4(q)$, where $q>2$. It is well known that $S$ includes a subgroup isomorphic to $G_2(q)$ (see, for example, \cite[Table 8.51]{13BHRD}). If $q$ is odd, then  $G_2(q)$ is unisingular
\cite[Theorem 1.3]{03GurTie}, and again using \cite{84Der}, we conclude that $p(q^2-1)\in\omega(H)\setminus\omega(S)$. Now let $q$ be even. We work similarly to the proof of Lemma 3.1 in \cite{13VasSt.t}. This proof involves three subgroups $A$, $B$ and $C$ of $G_2(q)$ such that  $B\simeq SL_3(q)$, $A$ is an SR-subgroup of type $A_1$ in $B$,  $C\simeq SL_2(q)$ and $C$ centralizes $A$. We choose $g\in A$ to be an element of order $q+1$. Since
 $g$ has a~nonzero fixed vector in any $B$-module over a field of characteristic 2 by Lemma \ref{l:sr}, it follows that $W=C_V(g)\neq 0$. The group $C$ acts on $W$ and includes a Frobenius subgroup
with cyclic kernel $F$ of order $q-1$ and complement of order 2. If $F\subseteq C_C(W)$, then $H$ has an element of order $2(q^2-1)$. Otherwise, we apply Lemma \ref{l:frob_action}
and derive that $H$ has an element of order $4(q+1)$. The centralizers of semisimple elements of $S$ are described in \cite{87DerMich}, and their structure implies that
$\omega_{mix}(S)$ consists of the divisors of the numbers $2(q^3+1)$, $2(q^3-1)$, $4(q^2+q+1)$, and $4(q^2-q+1)$. Since $q>2$, it follows that neither  $2(q^2-1)$ nor $4(q+1)$ divides any
of these numbers, and hence $\omega(H)\neq\omega(S)$.

Let $S=E_6^\pm(q)$ or $S=E_7(q)$. Setting $t=r_8(q)$ in the former case and $t=r_9(q)$ in the~latter, we see that $pt\in\omega(H)$ by Lemma \ref{l:spectrumofc} and $pt\not\in\omega(S)$ by \cite[Proposition 3.2]{05VasVd.t}.

We turn now to symplectic and orthogonal groups. To check whether or not a given integer belongs to the spectrum of $S$, we use \cite{10But.t}.

Let $n\neq 4$ and let $S$ be one of the groups $O_{2n+1}(q)$, with $q$ odd, or one the groups $O_{2n+2}^\pm(q)$. Define a number $t$ as follows.
If $n$ is even and $n\geqslant 8$, then take integers $n_1$ and $n_2$ such that $n=n_1+n_2$, $n_1>n_2\geqslant 3$ and $(n_1,n_2)=1$, and set $t=r_{n_1}(q)r_{n_2}(q)$. For $n=6$, set
$t=r_3(q)r_6(q)$ if $q\neq 2$ and $t=63$ otherwise. If $n$ is odd, then $t=r_n(q)$. It is easily seen that $t\in\omega(SL_n(q))$,  therefore, $pt\in\omega(H)$ by Lemma \ref{l:spectrumofc}.  On the other hand,
using \cite[Proposition 3.1]{05VasVd.t} when $t$ is a prime and \cite{10But.t} in other cases, we calculate that $pt\not\in\omega(S)$.

Let $S=O_9(q), O_{10}^\pm(q)$, where $q$ is odd. If $q\equiv \epsilon\pmod 4$, then $p(q^2+1)(q-\epsilon)\not\in\omega(S)$. On the other hand,  $(q^2+1)(q-\epsilon)\in\omega(SL_4^\epsilon(q))$, and so  $p(q^2+1)(q-\epsilon)\in\omega(H)$ by Lemma \ref{l:spectrumofc}.

Let $S=O^\pm_{10}(q)$, where $q$ is even. Since $S$ includes $Sp_8(q)$,  it includes a Frobenius subgroup with cyclic kernel of order $q^4+1$ and cyclic complement of order $8$ by Lemma \ref{l:frob_groups2}.
The action of $S$ on $V$ is faithful, and it follows by Lemma \ref{l:frob_action} that $H$ has an element of order $16$. It remains to note that $16\not\in\omega(S)$.

Let $S=S_{2n}(q)$, where $n\geqslant 2$ and $q=p^k$ is odd. Take $G=C_n(q)_u$ and identify $S$ with $G/Z(G)$. Denote the root system of $G$ by $\Phi$ and enumerate the fundamental roots as in Table~\ref{tab:z}.
Let $\alpha_0$ be the highest root of $\Phi$. The extended Dynkin diagram of $\Phi$ is as follows.
\begin{figure}[ht]
\centerline{\includegraphics[draft=false]{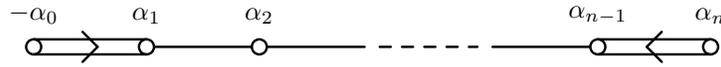}}
\vspace*{8pt}
\caption{Extended Dynkin diagram of type $C_n$}
\end{figure}

\noindent Observe that $\alpha_0=2(\alpha_1+\dots+\alpha_{n-1})+\alpha_n$, and using \cite[Theorem 1.12.1e]{98GorLySol}, we calculate that $h_{\alpha_0}(t)=h_{\alpha_1}(t)\dots h_{\alpha_n}(t)$.
Denote the subgroups generated by the root subgroups  associated with the roots $\pm\alpha_2,\dots,\pm\alpha_{n}$ and the roots $\pm \alpha_0$ by $A$ and $B$ respectively. Then $A$ and $B$ centralize each other, $A\simeq C_{n-1}(q)_u$ and $B\simeq A_1(q)_u$. Furthermore, a product of the central involutions $z_A$ and $z_B$ of $A$ and $B$, respectively,
is equal to the central involution $z$ of $G$. By Lemma \ref{l:sr}, every element of $A$ has a nonzero fixed vector in $V$. Choose $g\in A$ of order  $q^{n-1}-\epsilon$, where $\epsilon=1$ if $k(n-1)$ is odd and $\epsilon=-1$ otherwise, and set $W=C_V(g)$. The group $B$ centralizes $g$ and so acts on $W$. Observe that $z_B$ acting trivially because   $z_A$ is a power of $g$ and $z$ acts trivially on the whole $V$. Thus $W$ is a nonzero $B/Z(B)$-module.
Embedding $PSL_2(p)$ into $B/Z(B)$ and exploiting unisingularity of $PSL_2(p)$ \cite[Theorem 3.12]{01BabSha}, we
can take a element $h\in B$ of order $p+1$ such that $C_W(h)\neq 0$. Since $(q^{n-1}-\epsilon,p+1)=2$ by choice of
$\epsilon$, the image of $gh$ in $S$ has order $(q^{n-1}-\epsilon)(p+1)/2$, and hence $p(q^{n-1}-\epsilon)(p+1)/2\in\omega(H)$. It is not hard to see that $p(q^{n-1}-\epsilon)(p+1)/2\not\in\omega(S)$.

Let $S=S_{2n}(q)$, where $n\geqslant 4$ and $q$ is even. It is more convenient to consider $S$ as $B_n(q)_u$.
Again denote the root system of $S$ by $\Phi$,  enumerate the fundamental roots as in Table \ref{tab:z} and
denote the highest root of $\Phi$ by $\alpha_0$. 
\begin{figure}[ht]
\centerline{\includegraphics[draft=false]{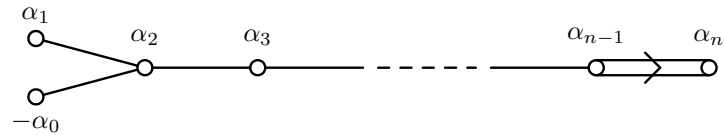}}
\vspace*{8pt}
\caption{Extended Dynkin diagram of type $B_n$}
\end{figure}

\noindent The set of long roots in $\Phi$ is a subsystem of type $D_n$, and the subgroup $A$ generated by the~corresponding root subgroups is isomorphic to $D_n(q)_u$.
Let $n$ be even. Denote the subgroups generated by the root subgroups associated with the roots $\pm \alpha_1,\pm\alpha_2,\dots,\pm\alpha_{n-2}$ and the roots $\pm \alpha_n$ by $B$ and $C$ respectively. It is clear that $B\simeq A_{n-2}(q)_u$, $C\simeq A_1(q)_u$ and $C$ centralizes $B$.  By Lemma \ref{l:sr}, every element of $B$ has a nonzero fixed vector in $V$. Choose $g\in B$ of order $r_{n-1}(q)$ and set $W=C_V(g)$. The group $C$ acts on $W$ and includes a Frobenius subgroup with cyclic kernel $F$ of order $q+1$ and complement of order $2$. If $F\subseteq C_C(W)$, then  $H$ has an element of order $2(q+1)r_{n-1}(q)$. Otherwise, Lemma \ref{l:frob_action} yields  $4r_{n-1}(q)\in\omega(H)$. 
None of the numbers $2(q+1)r_{n-1}(q)$ and $4r_{n-1}(q)$ lies in $\omega(S)$. Let $n$ be odd. Then take $B$ and $C$  to be the subgroups generated by the root subgroups associated with $\pm \alpha_1,\pm\alpha_2,\dots,\pm\alpha_{n-3}$ and $\pm \alpha_{n-1}, \pm \alpha_n$ respectively.
Then $B\simeq A_{n-3}(q)_u$, $C\simeq B_2(q)_u$ and again $C$ centralizes $B$. Choosing an element of order $r_{n-2}(q)$ in $B$ and working with a Frobenius group with cyclic kernel of order $q^2+1$ and cyclic complement of order $4$, which $C$ includes by Lemma \ref{l:frob_groups2}, we conclude that $H$ has either an element of order  $2(q^2+1)r_{n-2}(q)$ or an element of order $8r_{n-2}(q)$. And again none of these numbers lies in $\omega(S)$.

The case when $S=S_6(q)$ and $q$ is even requires some additional notation, so we extract this case as a separate lemma. Observe that this lemma is a part of Lemma 3.7 in  \cite{12Sta.t}, but
the~corresponding proof has a gap.

\setcounter{parag}{4}

\begin{lemma} \label{l:s6} Let $S=S_6(q)$, where $q$ is even. If $V$ is a nontrivial $S$-module over a field of characteristic $2$ and $H$ is a natural semidirect product of $V$ and $S$, then $\omega(S)\neq \omega(H)$.
\end{lemma}

\setcounter{parag}{1}

\begin{prf}
Let $K$ be the algebraic closure of the binary field and let $\overline G=Sp_6(K)$. Let $W$ be the natural $\overline G$-module and $f$ an nondegenerate  alternating bilinear form on $W$ preserved by $\overline G$.
There is a basis  $e_1$, \dots, $e_6$ of $W$ such that $f(e_i,e_j)=1$ if $i+j=7$ and $f(e_i,e_j)=0$ otherwise. Denote the set of matrices of $\overline G$ that are diagonal with respect to this basis by $\overline D$. Then
$\overline D$ is a~maximal torus of $\overline G$, and $\varepsilon_1$, $\varepsilon_2$, $\varepsilon_3$ defined by the formula $\varepsilon_i(d)=d(e_i)$ for $d\in \overline D$ are weights of $\overline G$ with respect to
this torus. The weights $\lambda_1=\varepsilon_1$, $\lambda_2=\varepsilon_1+\varepsilon_2$, and $\lambda_3=\varepsilon_1+\varepsilon_2+\varepsilon_3$ are the~fundamental weights of $\overline G$.
Given a dominant weight $\lambda$, we denote the irreducible $\overline G$-module with high weight $\lambda$ by $L(\lambda)$.

Identify $S$ with $C_{\overline G}(\sigma)$, where $\sigma$ is the Frobenius endomorphism of $\overline G$ that raises matrix entries to the $q$-th power. By \cite[Lemma 9]{08Zav1.t}, we may assume that the action of $S$ on $V$ is absolutely irreducible. By Steinberg's theorem \cite[13.3]{68Ste}, the module $V$ can be obtained as a~restriction of some irreducible $\overline G$-module $\overline V$ to $S$. Denote the high weight of $\overline V$ by $\lambda$. We show that $\omega(H)\setminus \omega(S)$ contains an element depending on $\lambda$. Note that the even numbers lying in $\omega(S)$ are exactly divisors of $8$, $4(q\pm 1)$, and $2(q^2\pm 1)$.

Suppose that $\lambda$ is neither $2^j\lambda_1$ nor $2^j\lambda_3$, with $j\geqslant 0$. Then by \cite[Theorem 1]{04Sup}, for every unipotent element $u$ of $\overline G$, the degree of the minimal polynomial
of $u$ on $\overline V$ is equal to $|u|$. Therefore, for every unipotent element $u$ of $S$,  the degree of the minimal polynomial
of $u$ on $V$ is equal to $|u|$ as well. Taking $u\in S$ to be a unipotent element of order $8$, we see that the~coset $Vu$ has an element of order 16.

Suppose that $\lambda=2^j\lambda_3$, where $j\geqslant 0$. 
Let $x\in S$ be an element of order $q^2+q+1$, and $\zeta\in K$ a primitive  $(q^2+q+1)$-th root of unity. Then $x$ is conjugate
to the diagonal matrix $y$ such that $\varepsilon_i(y)=\zeta^{q^{i-1}}$ for $1\leqslant i\leqslant 3$ in $\overline G$.  Since $\lambda(y)=\lambda_3(y)^{2^j}=(\zeta^{1+q+q^2})^{2^j}=1$, it follows that $y$
has a nonzero fixed vector in $\overline V$. Then $x$ also has a nonzero fixed vector in $\overline V$. This implies that $x$ has a nonzero fixed vector in $V$, and so the coset $Vx$ contains an element of order $2(q^2+q+1)$.

Lastly, let $\lambda=2^j\lambda_1$, where $j\geqslant 0$. Assume, for a while, that $j=0$, and hence $V$ is the natural module of $S$.
Exploiting the standard embedding of $Sp_2(q)\times Sp_4(q)$ into $S_6(q)$, we can find commuting elements $x\in S$ and $y\in S$ of orders 3 and $4$ respectively such that $\dim C_V(x)=4$ and the degree of the minimal polynomial of  $y$ on $C_V(x)$ is equal to 4. An easy computation establishes that the coset $Vxy$ contains an element of order $24$. Now let $j>0$. By Steinberg's tensor product theorem \cite{68Ste},
we have $L(\lambda)\simeq L(\lambda_1)^{(\rho^j)}$, where $\rho$ is the Frobenius endomorphism of $\overline G$ rasing  matrix entries to the second power. Since
twisting by $\rho$ does not affect the dimension of the fixed subspace or the degree of the minimal polynomial, the previous argument works, and again $H$ contains an element of order $24$.

Thus at least one of the numbers $16$, $2(q^2+q+1)$, and $24$ lies in $\omega(H)\setminus\omega(S)$, and so the~lemma follows.
\end{prf}

Finally, let $S=S_4(q)$, where $q>2$ is even. By Lemma \ref{l:frob_groups2}, there is a Frobenius subgroup with kernel of order $q^2+1$ and cyclic complement of order  $4$ in $S$. It follows that   $8\in\omega(H)\setminus\omega(S)$, and this completes the proof of Proposition \ref{p:equi}.

\end{prf}

\begin{prop}
Let  $S$ be one of the simple groups  $^3D_4(q)$, $F_4(q)$, $E_6(q)$, $^2E_6(q)$, $E_7(q)$, and $S_{4}(q)$. If $V$ is a nontrivial $S$-module over a field of characteristic $r$, where $r$ is coprime to $q$, and $H$ is a natural semidirect product of $V$ and $S$, then $\omega(S)\neq \omega(H)$.
\end{prop}

\begin{prf}

Similarly to the proof of Proposition \ref{p:equi}, we can assume that the action of $S$ on $V$ is faithful.
Denote the $r$-exponent of $S$ by $s$.

Let $S$ be one of the groups ${}^3D_4(q)$ and $F_4(q)$.
By \cite{81Wil, 89Kon.t}, the prime graph of $S$ is disconnected with $\pi(q^4-q^2+1)$ being a component. Suppose that $\omega(H)=\omega(S)$.
Then the prime graph of $H$ is also disconnected, and by \cite[Theorem A]{81Wil}, we see that $r$ and 2 lie in the same component of this graph. Since $q^4-q^2+1$ is odd, it follows that $r$ does not divide $q^4-q^2+1$
and $r(q^4-q^2+1)\not\in\omega(H)$. Exploiting the fact that $^3D_4(q)$ is a subgroup of $F_4(q)$ (see, for example, \cite{78Ste}) and applying \cite[Proposition 2]{13Zav.t},
we conclude that some element of $S$ of order $q^4-q^2+1$ has a nonzero fixed vector in $V$. Thus $r(q^4-q^2+1)\in\omega(H)$, and this is a contradiction.

Let $S=E^\varepsilon_6(q)$ and set $d=(3,q-\varepsilon)$. By \cite[Theorem 1]{13But.t}, the set $\omega_{p'}(S)$ consists of the~divisors of the numbers
\begin{multline}\label{eq:e6}\frac{q^6+\varepsilon q^3+1}{d},\frac{(q^4-q^2+1)(q^2+\varepsilon q+1)}{d}, \frac{(q^5-\varepsilon)(q+\varepsilon)}{d},%
q^5-\varepsilon, \\  \frac{(q^4+1)(q^2-1)}{d},  \frac{q^6-1}{d},  (q^3-\varepsilon)(q+\varepsilon), \frac{(q^4-1)(q^2-\varepsilon q+1)}{d}, q^4-1.\end{multline}

Let $\varepsilon=+$. The universal group $E_6(q)_{u}$ contains a subsystem subgroup $A$ isomorphic $A_5(q)_u$. By Lemma~\ref{l:frob_groups}, the group $A$
contains Frobenius subgroups of orders $q^5(q^5-1)_{(6,q-1)'}$ and $q^4(q^4-1)$. Clearly, these subgroups meet $Z(E_6(q)_{u})$ trivially, and so they can be embedded into $S$.
Then by Lemma \ref{l:frob_action}, there are elements of orders $r(q^5-1)_{(6,q-1)'}$ and $r(q^4-1)$ in $H$.
Suppose that both of these numbers lie in $\omega(S)$. Examining the numbers in (\ref{eq:e6}) and
using Zsigmondy's theorem, we calculate that $r$ divides $(q^5-1)(q+1)$ and $q^2-q+1$. Since  $((q^5-1)(q+1),
q^2-q+1)=(q+1,3)$, it follows that $r$ is equal to 3 and  it divides $q+1$. Then $s=(q^3+1)_3$. A maximal parabolic subgroup of $S$ with
Levi factor of type $D_5(q)$ has an~abelian unipotent radical (see, for example \cite[Theorem
1]{97Vas.t}) and $s\in\omega(D_5(q))$, therefore, by Lemma \ref{l:hh}, there is an element of order $3s$ in $H$. Thus $\omega(H)\neq\omega(S)$.

Let $\varepsilon=-$. Considering a subsystem subgroup $A_2(q^2)_u$ in $^2E_6(q)_u$ and applying Lemma \ref{l:frob_groups},
we conclude that $S$ has a Frobenius subgroup of order $q^4(q^4-1)_{(3,q^2-1)'}$. By Lemma \ref{l:frob_action} again,
we have $r(q^4-1)_{(3,q^2-1)'}\in\omega(H)$. Suppose that $r(q^4-1)_{(3,q^2-1)'}\in\omega(S)$. Working with the~2-parts of the numbers in (\ref{eq:e6}),
we deduce that $r$ divides $(3,q^2-1)(q^2+q+1)$, and in particular, $r$ is odd. Assume that $r$ divides $q^2+q+1$. Then
$s=(q^3-1)_r$. Consider a subsystem subgroup of $^2E_6(q)_u$ isomorphic to $^2A_5(q)_u$. Thus subgroup includes
$Z({}^2E_6(q)_{u})$ (see Table \ref{tab:z}), and hence its image in $S$ is an extension of the group of order $(2,q-1)$
by $^2A_5(q)$. By \cite[Lemma 5(4)]{11Gr}, the group $^2A_5(q)$ contains a Frobenius subgroup with kernel of order $q^3$ and cyclic complement of order $(q^3-1)_{(2,q-1)'}$.
It follows that $S$ also contains such a Frobenius subgroup, and so $r(q^3-1)_{2'}\in\omega(H)$.  Thus $rs\in\omega(H)\setminus\omega(S)$.
Now assume that $r$ does not divide $q^2+q+1$. Then $r=3$ and it divides $q^2-1$. If $q-1$ is divisible by $3$, then so is $q^2+q+1$, contrary to our assumption.
So 3 divides $q+1$ and $s=(q^3+1)_3$. A maximal parabolic subgroup of $S$ with Levi factor of type  $^2D_4(q)$ contains an element of order $s$, therefore, if $q$ is odd, then
we can apply Lemma \ref{l:hh} to deduce that $rs\in\omega(H)\setminus\omega(S)$. If $q$ is even, then  
$S$ contains a subgroup isomorphic to $^2D_4(q)$. A maximal parabolic subgroup of $^2D_4(q)$
with Levi factor of type $^2A_3(q)$ has an abelian unipotent radical (see, for example, \cite[Table 8.52]{13BHRD}) and contains an element of order $s$.
Applying Lemma \ref{l:hh}, we see that $rs\in\omega(H)\setminus\omega(S)$ in this case too.

Let $S=E_7(q)$. The structure of maximal tori of $E_7(q)_{u}$ described in \cite{91DerFak} implies that
 $\omega_{p'}(E_7(q)_{u})$ consists of the divisors of the numbers
\begin{multline}\label{eq:e7}(q^6+\epsilon q^3+1)(q-\epsilon), q^7-\epsilon, (q^4-q^2+1)(q^3-\epsilon), (q^5-\epsilon)(q^2+\epsilon q+1), \\%
(q^5-\epsilon)(q+\epsilon), \dfrac{(q^4+1)(q^2+1)(q-\epsilon)}{(2,q-1)},  (q^4+1)(q^2-1), (q^4-1)(q^2+\epsilon q+1), q^6-1,\end{multline}
where $\epsilon$ runs over $\{+1,-1\}$.
Since $E_7(q)$ contains a quotient of $SL_8(q)$ by a central $2$-subgroup, the group $S$ contains Frobenius subgroups of orders
$q^7(q^7-1)_{(2,q-1)'}$ and $q^4(q^4-1)$. It follows that $r(q^7-1)_{(2,q-1)'}$ and $r(q^4-1)$ lie in $\omega(H)$.
On the other hand, examining the numbers in  (\ref{eq:e7}), we see that
$r(q^4-1)\not\in\omega(S)$ if $r=2$ and $r(q^7-1)_{(2,q-1)'}\not\in\omega(S)$ if $r$ is odd.

Lastly, let $S=S_4(q)$, where $q>2$. By \cite[Lemma 2.8]{09HeShi}, there is a Frobenius subgroup with kernel of order $q^2$ and cyclic complement of order $(q^2-1)/(2,q-1)$ in $S$. By Lemma \ref{l:frob_action} yet again,
we have $r(q^2-1)/(2,q-1)\in\omega(H)$. By \cite{10But.t}, the number $(q^2-1)/(2,q-1)$ is maximal under divisibility in $\omega(S)$, and the proof of Proposition \ref{p:cross} is complete.
\end{prf}

Now we are ready to prove the theorem. Let $S$ and $H$ be as in the hypothesis and let $S$ be defined over a field of characteristic $p$. By Lemma \ref{l:reduction}, we can assume that $H$ is a semidirect product
of a nonzero finite $S$-module over a field of characteristic $r$ by $S$. If $r=p$, then $\omega(H)\neq\omega(S)$ by Proposition \ref{p:equi}.
Let $r\neq p$. If $S$ is an orthogonal or symplectic group other than $S_4(q)$, then $\omega(H)\neq\omega(S)$ by \cite[Theorem 1]{11Gr}. Otherwise, $\omega(H)\neq\omega(S)$
by Proposition \ref{p:cross}, and the theorem follows.

\setcounter{parag}{4}

\sect{Covers of $O_{2n+1}(q)$ and the spectrum of $S_{2n}(q)$}

A simple nonabelian group $L$ is said to be {\it quasirecognizable by spectrum} if every finite group $G$ with $\omega(G)=\omega(L)$ has exactly one nonabelian composition factor and this factor is isomorphic to $L$.
The following result was obtained in studying quasirecognizability of simple symplectic groups: If $L=S_{2n}(q)$, where $n\geqslant 4$, and $G$ is a finite group such that $\omega(G)=\omega(L)$, then $G$ has exactly
one nonabelian composition factor and if this factor is a group of Lie type in the same characteristic as $L$, then it is isomorphic either to $L$ or to one of the groups $O_{2n+1}(q)$ and $O_{2n}^-(q)$ \cite[Theorem 3]{09VasGrMaz.t}.
As we mentioned in Introduction, the case when the factor is isomorphic to $O_{2n+1}(q)$ requires a separate treatment because the~spectra of $S_{2n}(q)$ and $O_{2n+1}(q)$ are very close.

\setcounter{parag}{1}
\setcounter{prop}{2}

\begin{prop}
Let $L=S_{2n}(q)$ and $S=O_{2n+1}(q)$, where $n\geqslant 3$ and $q$ is odd. If $H$ is a proper cover of $S$, then $\omega(H)\not\subseteq \omega(L)$.
In particular, if $G$ is a finite group such that $\omega(G)=\omega(L)$ and $G$ has a composition factor isomorphic to $S$, then $S\leqslant G\leqslant \operatorname{Aut}(S)$.
\end{prop}

\setcounter{parag}{5}

\begin{prf}
Let $q$ be a power of a prime $p$. By Lemma \ref{l:reduction}, we can assume that $H=V\leftthreetimes S$, where $V$ is a nonzero  $S$-module over a field of a positive characteristic $r$. Suppose that $r\neq p$ and $(S,r)\neq (O_7(3),5)$.
Then by the proof of Proposition 4 in \cite{11Gr}, the group $H$ contains an element of order $rt$ such that $p$ is coprime to $t$ and $rt\not\in\omega(S)$. Since
the isomorphism types of maximal tori in $S$ are the same as in $L$ \cite{81Car}, we have $\omega_{p'}(S)=\omega_{p'}(L)$, and so $rt\not\in\omega(L)$. If $S=O_7(3)$ and $r=5$, then
the~5\nobreakdash-Brauer character table of $S$ \cite{95AtlasBr} implies that every element of $S$ of order 7 has a~nonzero fixed vector in $V$, and hence $35\in\omega(H)\setminus\omega(L)$.
Now let $r=p$. If $n\neq 4$, then we showed in the proof of Proposition \ref{p:equi} that $H$ contains an element of order $pt$, where $t$ is odd and $pt\not\in\omega(S)$.
One can show (by using, for example, \cite{10But.t}) that the subsets of odd numbers are the same in $\omega(S)$ and $\omega(L)$, and thus $pt\not\in\omega(L)$.
If $n=4$ and $q\equiv \epsilon\pmod 4$, then we proved that $p(q^2+1)(q-\epsilon)\in\omega(H)$. By \cite{10But.t}, we see that $p(q^2+1)(q-\epsilon)\not\in\omega(L)$.

Now suppose that $G$ is a finite group such that $\omega(G)=\omega(L)$. By  \cite[Corollary 7.2]{05VasVd.t}, the~group $G$ has exactly one nonabelian composition factor and this factor is isomorphic to $S$.
Denoting the soluble radical of $G$ by $K$, we see that $S\leqslant G/K\leqslant \operatorname{Aut}(S)$. If $K\neq 1$, then $G$ contains a proper cover of $S$ and so $\omega(G)\not\subseteq \omega(L)$, which contradicts the hypothesis.
Thus $K=1$, and the proof is complete.
\end{prf}

\end{document}